\magnification=\magstep1
\vbadness=10000
\hbadness=10000
\tolerance=10000

\def\F{{\bf F}}  
\def\m{{\bf m}}  
\def\Q{{\bf Q}}  
\def\Z{{\bf Z}}  

\proclaim The fake monster formal group. \hfill 27 May, 31 Dec 1998

Richard E. Borcherds, 
\footnote{$^*$}{ Supported by a Royal Society
professorship and by an NSF grant.}

D.P.M.M.S.,
16 Mill Lane,
Cambridge,
CB2 1SB,
England.

e-mail: reb@dpmms.cam.ac.uk

home page: www.dpmms.cam.ac.uk/\~{}reb

\bigskip

\proclaim Contents.

1. Introduction.

Notation and terminology.

2.~Some theorems about smooth Hopf algebras.

3.~Liftings of Lie algebra elements.

4.~The smooth fake monster Hopf algebra.

5.~A smooth Hopf algebra for the Virasoro algebra.

6.~The no-ghost theorem over $\Z$.

7.~An application to modular moonshine.

8.~Open problems.

\proclaim 1.~Introduction.

The main result of this paper is the construction of
``good'' integral forms for the
universal enveloping algebras of
the fake monster Lie algebra and the Virasoro
algebra. As an application we construct formal group laws over
the integers for these Lie algebras. We also prove a form of the
no-ghost theorem over the integers, and use this to verify an
assumption used in the proof of the modular moonshine conjectures.

Over the integers the universal enveloping algebra of a Lie algebra is
not very well behaved, and it is necessary to use a better integral
form of the universal enveloping algebra over the rational
numbers. The correct notion of ``good'' integral form was found by
Kostant [K].  He found that the good integral forms are the ones with
a structural base (as defined in 2.3), and showed that finite
dimensional semisimple Lie algebras have a structural base.  The
existence of a structural base implies that the dual algebra of the
underlying coalgebra is a ring of formal power series. As this ring
can be thought of as a sort of ``coordinate ring'' of some sort of
formal group, this condition can be though of as saying that the
formal group is smooth and connected, and that the formal
group comes from a ``formal group law''.

The main point of this paper is to find such integral forms for
the universal enveloping algebras of certain infinite dimensional 
Lie algebras. We recall that the universal enveloping algebra of a Lie algebra
has a natural structure of a cocommutative Hopf algebra. 
So we need some theorems to tell us when a Hopf algebra has a structural
base.  In section 2 we prove that a Hopf algebra $H$ has a structural
base provided that there are ``sufficiently many'' group-like elements
$1+a_1x+a_2x^2+\cdots\in H[[x]]$, and provided a few other minor
conditions are satisfied.  Such a group-like element should be thought
of as an infinitesimal curve in the formal group, so roughly speaking
this condition means that a set of generators of the Lie algebra of
the formal group can be lifted to formal curves in the formal
group. As an example of the theorems in section 2, we find that the
Lie bracket $[a_1,b_1]$ of two primitive liftable elements is also
liftable. At first sight this looks as if it should be easy to prove,
as all we need to do is write down an explicit lifting, with
coefficients that are some universal (non-commutative) polynomials in
the coefficients of liftings of $a_1$ and $b_1$. However this seems to
be very hard to do explicitly (partly because such a lifting is far
from unique, which paradoxically makes it harder to find one). Instead,
we use a far more roundabout argument, where we first need to
prove a theorem (theorem 2.12) saying that all primitive elements of
certain Hopf algebras are liftable.

We also give an incidental application of these theorems to certain Hopf
algebras $F_n$ considered by Dieudonn\'e. He showed that these Hopf
algebras have a structural basis when reduced mod $p$ for any prime
$p$; we show that the Hopf algebras $F_n$ already have a structural
basis over the integers.

In section 2 we have reduced the problem of finding a structural basis 
of a Hopf algebra $H$
to the problem of finding enough group-like elements in $H[[x]]$. 
For the universal enveloping algebras of 
finite dimensional simple Lie algebras, or more generally
for Kac-Moody algebras, this is easy to do, because
we just take the formal 1-parameter subgroup
$$\sum_{n\ge 0} x^n Ad(e)^n/n!$$ for (locally nilpotent) elements $e$
of the Lie algebra corresponding to real roots, and define
the integral form of the Hopf algebra to be the one generated by the coefficients of these liftings. This gives the usual
Kostant integral form for the universal enveloping algebra 
of finite dimensional semisimple Lie algebras.
For the Lie algebras in this paper this does not work because there
are not enough locally nilpotent elements. (We can try to use the
elements $e^n/n!$ for non nilpotent elements, but there seems no
obvious way in which this gives a good integral form.)  Fortunately we
do not need to lift the generators of the Lie algebra to formal one
parameter groups, and it is sufficient to lift them to formal
curves, which is easier to do.

In section 3 we construct liftings of some elements of the vertex
algebra of a lattice to formal curves. More precisely, we show that we
can lift any element of a certain Lie algebra provided that it lies in
a root space of a root of norm 2 or 0. We do this by explicitly writing
down a group-like lifting, and checking by brute force that its
coefficients are integral. (Unfortunately there seems no obvious way
to extend this brute force approach to root spaces of negative norm
roots.)

In section 4 we use the liftings of section 3 to construct a smooth
integral form of the universal enveloping algebra of the fake monster
Lie algebra, which is constructed from the lattice $II_{25,1}$. The
main point about this Lie algebra is that it is generated by the root
spaces of norm 2 and norm 0 roots, so in section 3 we have constructed
enough liftings to apply the theorems in section 2. (This is a very
special property of $II_{25,1}$; it is the only known indefinite
lattice whose Lie algebra is generated by the root spaces of norm 2
and norm 0 roots. In other words the theory in sections 2 and 3 has
been developed mainly for this one example!)  We can summarize the
main results about the fake monster Lie algebra proved in this paper
as follows.
\proclaim Theorem 4.1. There is a $II_{25,1}$-graded
Hopf algebra $U^+(\m)$ over $\Z$ with the following
properties.
\item {1} $U^+(\m)$ has a structural basis over $\Z$.
\item {2} The primitive elements of $U^+(\m)$ are an integral form of the
fake monster Lie algebra $\m$.
\item {3} For every norm 2 vector of $II_{25,1}$, $U^+(\m)$ contains
the usual (Kostant) integral form of the universal enveloping algebra
of the corresponding $sl_2(\Z)$.

In section 5 we construct a smooth integral form for the universal
enveloping algebra of the Virasoro algebra, by applying the theorems
in section 2 to a set of explicit liftings of the basis elements of
the Virasoro algebra. In other words we construct a formal group law
for the Virasoro algebra over the integers. 

The no-ghost theorem in string theory states that a certain real
vector space of states is positive definite (so it contains no
negative norm vectors, which are sometimes called ghosts and which
would prevent the space from being a Hilbert space).  This real vector
space has a natural integral form which can be made into a positive
definite lattice using the inner product, and we can ask for an
integral form of the no-ghost theorem, which should say something
about the structure of this lattice.  In section 6 we use the smooth
formal group of the Virasoro algebra to prove an integral form of the
no-ghost theorem, at least in the case of vertex algebras constructed
from 26 dimensional lattices. More precisely we prove bounds on the
determinant of this lattice, which quite often imply that the lattice
is self dual.

The proof [B98] of Ryba's modular moonshine conjectures for large
primes used an unproved technical assumption about the monster Lie
algebra.  In section 7 we use the integral no-ghost theorem to prove
this assumption, thus completing the proof of the modular moonshine
conjectures for odd primes. (The proof in [B-R] for the prime 2 still
relies on another unproved technical assumption.)

\proclaim Notation and terminology.

\item{$\alpha$} An element of a lattice.
\item{$\beta$} An element of a lattice.
\item{$\gamma$} A norm 0 element of a lattice.
\item{$\Gamma$} An indeterminate.
\item{$c$} An element generating the center of the Virasoro algebra. 
\item{$Der$} The derivations of a commutative ring. 
\item{$e^\alpha$} An element of the group ring of a lattice $L$, usually
regarded as an element of the vertex algebra of $L$. 
\item{$\epsilon_i$} The function from $I$ to $\Z$ taking value 1 on $i\in I$
and 0 elsewhere.
\item{$\F_p$} A finite field of order $p$.
\item{$F(\lambda)$} If $\lambda=1^{i_1}2^{i_2}\cdots$ is a partition,
then $F(\lambda)=i_1!i_2!\cdots$.
\item{$\gamma$} A norm 0 vector of a lattice $L$.
\item{$H$} A cocommutative Hopf algebra.
\item{$I$} An index set, or a finite sequence of integers.
\item{$II_{m,n}$} An even self dual lattice of dimension  $m+n$
and signature $m-n$.
\item{$J$} A finite sequence of integers.
\item{$l(I)$} The length of the sequence $I$. 
\item{$l(\lambda)$} If $\lambda=1^{i_1}2^{i_2}\cdots$ is a partition,
then $l(\lambda)=i_1+i_2+\cdots$.
\item{$\lambda$} A partition.
\item{$K$} A field, sometimes the quotient field of $R$.
\item{$L$} A lattice, usually even and self dual.
\item{$L_i$} A basis element for the Virasoro or Witt algebra. 
\item{$\m$} An integral form of the fake monster Lie algebra
(or the monster Lie algebra in section 7).
\item{$N$} An integer.
\item{$P(\lambda)$} If $\lambda=1^{i_1}2^{i_2}\cdots$ is a partition,
then $P(\lambda)$ is the integer $1^{i_1}2^{i_2}\cdots$.
\item{$p$} A prime.
\item{$p(n)$} The number of partitions of $n$. 
\item{$\Q$} The rational numbers.
\item{$\rho$} A Weyl vector. 
\item{$R$} A commutative ring.
\item{$S$} The antipode of a Hopf algebra.
\item{$\Sigma(I)$} The sum of the elements of $I$. 
\item{$U$} A universal enveloping algebra.
\item{$U^+$} An integral form of a universal enveloping algebra
with a structural basis.
\item{$V$} A vertex algebra.
\item{$V_p$} The shift (or Verschiebung) of a coalgebra over $\F_p$.
\item{$Vir$} The Virasoro algebra.
\item{$Witt$} The Witt algebra.
\item{$x$} A formal variable.
\item{$\Z$} The integers. $\Z_{\ge0}$ is the set of non-negative integers.
\item{$\Z_p$} The $p$-adic integers. 
\item{$Z_\alpha$} An element of a structural basis.
\item{$z$} A formal variable.
\bigskip
\item{}{\bf Bialgebra} A module with compatible algebra and coalgebra
structures.
\item{}{\bf Group-like} $\Delta(a)=a\otimes a$.
\item{}{\bf Hopf algebra} A bialgebra with antipode.
\item{}{\bf Liftable} See definition 2.2.
\item{}{\bf Primitive} $\Delta(a) = a\otimes 1+1\otimes a$.
\item{}{\bf Smooth.} A Hopf algebra or coalgebra with a structural basis.
\item{}{\bf Structural basis} $\Delta Z_\alpha=\sum_{0\le \beta \le \alpha}
Z_\beta\otimes Z_{\alpha-\beta} $

\proclaim 2.~Some theorems about Hopf algebras.

In this section we will prove several results about Hopf algebras. The
main result is theorem 2.15, which states that under mild conditions a
bialgebra generated by the coefficients of group-like liftings has
a structural basis. In later
sections we will use this to construct a good integral form of the
universal enveloping algebra of the fake monster Lie algebra.

We recall that a coalgebra over a commutative ring $R$ is
an $R$-module $H$ with a coassociative coproduct and a counit,
a bialgebra is a module with compatible algebra and coalgebra structures,
and a Hopf algebra is a bialgebra with an antipode $S$.

\proclaim Definition 2.1. An element $a$ of a coalgebra is called 
group-like if $\Delta(a)=a\otimes a$, and primitive if
$\Delta(a)=a\otimes 1+1\otimes a$.

\proclaim Definition 2.2. We say that an element $a_1$ of
the coalgebra $H$ is liftable to a group-like element, or liftable
for short, if we can find elements $a_n\in H$ $(n\ge 0)$ such that $a_0=1$
and the element $\sum_{n\ge 0} a_nz^n\in H[[z]]$ is
group-like. In other words
$$\Delta(a_n)=\sum_{0\le m\le n} a_m\otimes a_{n-m}$$
for all $n\ge 0$.
We say that $a_1$ is liftable to order $N$ if we can find elements
$a_n$ for $0\le n\le N$ (with $a_0=1$) satisfying the relations above
for $0\le n\le N$. We say that the lifting is graded if $H$
is graded by some abelian group $L$ and $\deg(a_n)=n\alpha$
for some $\alpha\in L$.

We write $\Z_{\ge 0}$ for the nonnegative integers, and $\Z_{\ge
0}^{(I)}$ for the functions from a set $I$ to $\Z_{\ge 0}$ that are
zero on all but a finite number of elements of $I$. The element
$\epsilon_i$ of $\Z_{\ge 0}^{(I)}$ is defined to be the function that
is 1 on $i\in I$ and 0 elsewhere. We define a partial order $\ge$ on
$\Z_{\ge 0}^{(I)}$ in the obvious way, by saying $\alpha\ge \beta$ if
$\alpha(i)\ge \beta(i)$ for all $i\in I$.

\proclaim Definition 2.3. A structural basis for a coalgebra or
Hopf algebra
$H$ over a ring $R$ is a set of elements $Z_\alpha$,
$\alpha\in \Z_{\ge 0}^{(I)}$ for some set $I$,
such that the elements $Z_\alpha$ form a basis for the free $R$-module $H$
and such that $\sum_\alpha Z_\alpha x^\alpha$ is group-like, in other words
$$\Delta(Z_\alpha) =
\sum_{0\le \beta\le \alpha}Z_\beta\otimes Z_{\alpha-\beta}.$$

\proclaim Lemma 2.4. If a coalgebra $H$ has a structural basis $Z_\alpha$
then the $R$-module of primitive elements of $H$ has a basis consisting of
the elements $Z_{\epsilon_i} $ for $i\in I$.

Proof. We have to show that every primitive element of $H$ is a linear
combination of the elements $Z_{\epsilon_i}$. Give $H$ the bialgebra
structure such that $Z_\alpha Z_\beta=
Z_{\alpha+\beta}
\prod_{i\in I}{\alpha(i)+\beta(i)\choose\beta(i)}
$ (as in [A,
section 2.5.1] when $R$ is a field). The dual algebra $H^*$ is a ring
of formal power series and $H$ acts on $H^*$ as differential
operators. Suppose $D$ is any primitive element of $H$. By
subtracting multiples of $Z_{\epsilon_i}$ from $D$ we can assume that
$D$ acts trivially on the elements $x_i$ dual to $Z_{\epsilon_i}$.
But then the fact that $D$ is a derivation implies that $D$ acts
trivially on any polynomial in the $x_i$'s. As there are no elements
of $H$ orthogonal to all polynomials in the $x_i$'s this proves that
$D$ must be 0. This proves lemma 2.4.

\proclaim Lemma 2.5. The liftable elements of a bialgebra over $R$ form
an $R$-submodule.

Proof. If $a(x)$ and $b(x)$ are lifts of $a_1$ and $b_1$
then $a(x)b(x)$ is a lift of $a_1+b_1$, and if $r\in R$ then
$a(rx)$ is a lift
of $ra_1$. This proves lemma 2.5.

\proclaim Lemma 2.6. If $H$ is a bialgebra with a structural
basis then all primitive elements
of $H$ are liftable.

Proof. By lemma 2.4 the $R$-module of primitive elements has a base of
elements $Z_{\epsilon_i}$. Each of these basis elements can be lifted
by $\sum_{n\ge 0} Z_{n\epsilon_i}x^n$. Lemma 2.6 now follows from
lemma 2.5.

\proclaim Lemma 2.7. Suppose $H$ is a coalgebra with a structural basis
over an integral domain
$R$. Then any nonzero coideal $J$ of $H$ contains
a nonzero primitive element of $H$.

Proof.
We first assume that $R$ is a
field. A coalgebra over a field with a structural basis is
irreducible and pointed (by [A, section 2.5]), so lemma 2.7 follows
immediately from [A, corollary 2.4.14], which states that any nonzero
coideal of a pointed irreducible coalgebra over a field contains a
nonzero primitive element.

For the general case, let $K$ be the quotient field of $R$. If the
coideal $J$ is nonzero then $J\otimes_R K$ is a nonzero coideal of
$H\otimes_R K$. As lemma 2.7 is is true over the field $K$ this shows
that $J\otimes_R K$ contains a primitive element, and multiplying this
by a suitable nonzero constant to cancel the denominators gives a
nonzero primitive element of $J$.  This proves lemma 2.7.

\proclaim Lemma 2.8.
Suppose that $H$ is a torsion-free coalgebra over $\Z$ such that
for every prime $p$ every primitive
element of $H/pH$ is the image of a primitive
element of $H$. Let $N$ be any positive integer. Then every
primitive element of $H/NH$ is the image of a primitive element of
$H$.

Proof. We prove this by induction on the number of prime factors of
$N$, the case $N=1$ being trivial. Suppose $p$ is any prime factor of
$N$, and suppose that the image of $x\in H$ is primitive in $H/NH$,
so that we have to show that $x$ is congruent to a primitive element
mod $N$. The
element $x$ maps to a primitive element in $H/pH$, so by assumption we
have $x=y+pz$ for some primitive element $y$ and some $z$. But then
$\Delta(pz)\equiv pz\otimes 1+1\otimes pz\bmod N$ as $x$ and $y$
are both primitive mod $N$. Using the fact
that $H$ is torsion-free we can divide
through by $p$ to find that $\Delta(z)\equiv z\otimes 1+1\otimes
z\bmod N/p$. By induction on the number of prime factors of $N$
this shows that $z=t+(N/p)u$ for some
primitive element $t$. Substituting this back in to $x=y+pz$
shows that $x=y+pt+Nu$, so $x$ is congruent to a
primitive element $y+pt$ mod $N$. This proves lemma 2.8.

\proclaim Lemma 2.9. If every primitive element of a bialgebra
is liftable to all
finite orders, then every primitive element
is liftable (to infinite order).

Note that we do not claim that if one particular primitive element of
a cocommutative Hopf algebra can be lifted to all finite orders then
it can be lifted to infinite order. I do not know whether or not
this is true in
general.

Proof. Suppose that $1+a_1x+a_2x^2+\cdots + a_nx^n$ is an order $n$
lifting of $a_1$. It is sufficient to show that we can extend this to
an order $n+1$ lifting, because by repeating this
an infinite number of times we get a lifting of $a_1$ of infinite order.
We will prove by induction on $k$ that if $k\le
n$ then we can find an order $n+1$ lifting of $a_1$ which agrees with
$a(x)$ to order $k$. This is true for $k=1$ by the assumption that
$a_1$ has liftings of arbitrary finite order. Let $b(x)=1+b_1x+\cdots
+b_{n+1}x^{n+1}\in H[x]/(x^{n+2})$ be an order $n+1$
lifting of $a_1$ with $b_i=a_i$ for $i\le k$. Then
$a(x)b(x)^{-1}= 1+(a_{k+1}-b_{k+1})x^{k+1}+\cdots$ is group-like, so
$a_{k+1}-b_{k+1}$ is primitive. (Note that $b(x)^{-1}$ is a well
defined order $n+1$ lifting, as we can construct it as
$\sum_m (1-b(x))^m$.) Let $c(x)$ be an order $n+1$ lifting
of $a_{k+1}-b_{k+1}$ (which exists as we assumed all primitive
elements have liftings of any finite order). Then $b(x)c(x^{k+1})$ is
an order $n+1$ lifting of $a_1$, which agrees with $a(x)$ to order
$k+1$. This shows that any order $n$ lifting of $a_1$ extends to an
order $n+1$ lifting. This proves lemma 2.9.

\proclaim Lemma 2.10. Suppose that $H$ is a torsion-free bialgebra over $\Z$
and suppose that for every prime $p$ the map from
primitive elements of $H$ to primitive elements of $H/pH$ is onto.
Then every primitive element of $H$ is liftable.

Proof. We will show by induction on $n$ that any primitive element
$a$ can be lifted to a group-like element in $H[x]/(x^{n+1})$; this
is trivial for $n=0$.
We will show by induction
on $k$ that if $1\le k\le n$ then
we can find a lifting $\sum_{0\le i\le n}c_ix^i$
of $n!a$
such that $n!^i|c_i$ for $i\le k$. This is true for $k=1$,
because we can use the order $n$ lifting
$\sum_{0\le i\le n} (n!^i/i!)a^ix^i$.
So suppose that $k\ge 2$ and that $n!^i|c_i$ for $i<k$. Then
$$\Delta(c_k)=\sum_{0\le i\le k}c_i\otimes c_{k-i}
\equiv c_k\otimes 1+1\otimes c_k\bmod n!^k.$$
Therefore $c_k$ is primitive in $H/n!^kH$.
By lemma 2.8 and the assumption that for every prime $p$ the map from
primitive elements of $H$ to primitive elements of $H/pH$ is onto
this implies that $c_k=d_1+n!^ke$ for some primitive $d_1$
and some $e$. By induction on $n$ we can find a lifting
$d(x)=\sum_{0\le i<n} d_ix^i\in H[x]/(x^n)$ of $d_1$, and changing
$x$ to $x^k$ we get a group-like element $d(x^k)\in H[x]/(x^{n+1})$.
(Here we use the fact that $k\ge 2 $ and $n\ge 1$, so that $kn\ge n+1$.)
Multiplying $\sum c_ix^i$ by $d(x)^{-1}\bmod x^{n+1}$ we get a group-like
element in $H[x]/(x^{n+1})$ whose coefficient of $x^i$ for $i\le k$
is divisible by $n!^k$. This completes the proof by induction on $k$,
and shows that we can assume that $n!^i|c_i$ for all $i\le n$. But then
using the fact that $H$ is torsion free we see that
$\sum (c_i/n!^i)x^i$ is a lifting of $c_1/n!=a$ in $H[x]/(x^{n+1})$.

We have shown that any primitive element of $H$ can be lifted to a
group-like element of any finite order $n$. By lemma 2.9 this implies
that every primitive element of $H$ is liftable.  This proves lemma
2.10.

We recall that for any cocommutative
coalgebra or Hopf algebra $H$ over $\F_p$ there is a homomorphism
$V_p$ from $H$ to itself, called the shift (or
Verschiebung), which is dual to the Frobenius map
$x\mapsto x^p$ on the dual algebra $H^*$ (see for example [A section 5.4]).
The shift $V_p$ commutes with all homomorphisms of coalgebras.
If $H$ is a coalgebra over $\Z$ we also write $V_p$ for the shift
of $H/pH$.

\proclaim Lemma 2.11. Suppose that the $Z_\alpha$'s for
$\alpha\in \Z_{\ge 0}^{(I)}$ are elements
of a coalgebra $H$ over $\F_p$ satisfying the relations
$$\Delta(Z_\alpha)=\sum_{0\le \beta\le \alpha}Z_\beta\otimes
Z_{\alpha-\beta}$$ Then $V_p(Z_\alpha)=Z_{\alpha/p}$, where
$Z_{\alpha/p}$ means 0 if $p$ does not divided $\alpha$.

Proof. If the $Z_\alpha$'s are linearly independent and span $H$,
or in other words if they form a structural basis,
then the result
follows from [A, theorem 2.5.9]. In general we
can find a coalgebra with a structural basis
together with a homomorphism mapping the elements of
the structural basis to the $Z_\alpha$'s. The result then follows
because the shifts $V_p$ commute with homomorphisms. This proves lemma
2.11.

\proclaim Theorem 2.12. Suppose that $H$ is a $\Z_{\ge 0}$-graded
cocommutative
bialgebra over $\Z$ whose homogeneous pieces are finitely generated free
abelian groups
and such that the degree 0 piece is spanned by 1. Suppose
that the shifts $V_p:H/pH\rightarrow H/pH$
are onto for all primes $p$. Then every primitive element of $H$
is liftable.

Remark. In fact theorem 2.15 (whose proof uses theorem 2.12) implies the
stronger result that $H$ has a structural basis.

Proof. We first show that if $K$ is any field then the $K$ bialgebra
$C=H\otimes K$ is irreducible, in other words it has only one simple
subcoalgebra.  By theorem 2.3.4 of [A], the simple subcoalgebras of a
$K$-coalgebra $C$ correspond to the maximal ideals of the dual algebra
$C^*$ which are not dense, so it is sufficient to show that $C^*$ is a
local ring.  The elements of $C^*$ are infinite sums of the form
$c=c_0+c_1+\cdots$ with $c_i\in C_i^*$ (where $C_i$ is the degree $i$
subspace of $C$). We show that the ideal $M$ of all elements with
$c_0=0$ is the unique maximal ideal of $C^*$, in other words that any
element $c$ with $c_0\ne 0$ invertible. The space $C_0^*$ is
isomorphic to the field $K$, so if $c_0\ne 0$ then it is invertible,
so we can recursively define an inverse $d=d_0+d_1+\cdots$ to $c$ by
putting $d_0=c_0^{-1}$, $d_n=-c_0^{-1}(d_{n-1}c_1+\cdots+d_0c_n)$ for
$n>0$.  This shows that $C^*$ is a local ring, so $C$ is irreducible.

Next we check that $H\otimes \F_p$ and $H\otimes \Q$ have structural bases.
By [A, theorem 2.4.24] any irreducible bialgebra over a field,
in particular $H\otimes \F_p$,
is automatically a Hopf algebra.
By [A,corollary 2.5.15] any cocommutative irreducible
Hopf algebra over $\F_p$
such that $V_p$ is onto, in particular $H\otimes \F_p$,
has a structural basis. Similarly
theorem 2.5.3 of [A] states that any irreducible cocommutative Hopf
algebra over a field of characteristic 0 is a universal enveloping
algebra of its primitive elements. So $H\otimes \Q$ is the universal
enveloping algebra of its primitive elements and therefore also has a
structural basis.

The next step is to show that any primitive element of $H\otimes \F_p$
can be lifted to a primitive element of $H$.
If $m_\alpha$ is the dimension of the space of primitive elements
of degree $\alpha$
of a $\Z_{\ge 0}^n$-graded Hopf algebra with a structural basis
and $n_\alpha$ is the dimension
of the space of all elements of degree $\alpha$, then
$$\prod_{\alpha> 0} (1-x^\alpha)^{-m_\alpha}=\sum_{\alpha\ge 0} n_\alpha
x^\alpha$$
by lemma 2.4. In particular the spaces of
primitive elements of degrees $\alpha$ in $H\otimes \F_p$ and
$H\otimes \Q$ have the same dimension, because it is obvious that the
spaces of all elements in $H\otimes \F_p$ and $H\otimes \Q$ of degree
$\alpha$ have the same dimension, and $H\otimes \F_p$ and
$H\otimes \Q$ both have structural bases by the paragraph above.
The space $P_\alpha$ of primitive
elements of $H$ of degree $\alpha$ also has the same rank as the
space of primitive elements of $H\otimes \Q$, and therefore the same
rank as the space $P_{\alpha,p}$ of primitive elements of $H\otimes
\F_p$ of degree $\alpha$. The map from $P_\alpha/pP_\alpha$ to
$P_{\alpha,p}$ is obviously injective and both spaces are finite
dimensional vector spaces over $\F_p$ of the same dimension, so the
map between them is also surjective, or in other words every primitive
element of $H\otimes \F_p$ (of any degree $\alpha$) can be lifted to a
primitive element of $H$.
By lemma 2.10 this implies every primitive element of $H$
is liftable. This proves theorem 2.12.

\proclaim Corollary 2.13. Let $k$ and $n$ be positive integers.
If
$$a(x)=\sum_{i_1,\ldots, i_k\ge 0} a_{i_1,\ldots, i_k}x_1^{i_1}\cdots
x_k^{i_k}$$ is group-like, for $a_{i_1,\ldots, i_k}$ in some bialgebra
over a commutative ring $R$, and $a_{0,\ldots,0}=1$,
$a_{i_1,\ldots,i_k}=0$ for $0<i_1+\cdots+i_k< n$, then
$a_{i_1,\ldots,i_k}$ is liftable for $i_1+\cdots+i_k=n$.

Proof. It is sufficient to prove this for the universal example, where
the bialgebra $H$ is the free associative algebra over $\Z$ generated
by $a_{i_1,\ldots,i_k}$ for $i_1+\cdots,i_k\ge n$ and the coalgebra
structure is defined so that $a(x)$ is group-like. Put $a_0=1$,
$a_{i_1,\ldots,i_k}=0$ for $0<i_1+\cdots i_k<n$, and grade $H$ by
letting $a_{i_1,\ldots,i_k}$ have degree $i_1+\cdots+i_k$. Then
$V_p(a_{i_1.\ldots,i_k})=a_{i_1/p,\ldots, i_k/p}$ by lemma 2.11, so
$V_p:H/pH\mapsto H/pH$ is onto. By theorem 2.12 this implies that any
primitive element of $H$, in particular $a_{i_1,\ldots,i_k}$ for
$i_1+\cdots +i_k<2n$, is liftable. This proves corollary 2.13.

\proclaim Theorem 2.14. The set of liftable elements of a cocommutative
bialgebra $H$ over a commutative ring $R$ is a Lie algebra
over $R$.

Proof. We already know that the liftable primitive elements form a
module over $R$ by lemma 2.5.
We have to show that the set of liftable elements is
closed under the Lie bracket. Suppose $a_1$ and $b_1$ are two
primitive elements with group-like lifts $a(x)=\sum a_ix^i$ and
$b(y)=\sum_ib_iy^i$. It is sufficient to prove theorem 2.14 in the
case of the universal example, when $H$ is the free associative
algebra over $\Z$ generated by the elements $a_i$ and $b_i$ for $i>0$,
with the coalgebra structure defined by the fact that $a(x)$ and
$b(y)$ are group-like. We $\Z$-grade $H$ by giving $a_i$ and $b_i$
degree $i\in \Z$. By lemma 2.11 the shifts $V_p$ contain all the
elements $a_i$ and $b_i$ in their images, so the shifts are onto as
they are homomorphisms of $\F_p$-algebras. Then by theorem 2.12 we see
that all primitive elements of $H$, in particular $[a_1,b_1]$,
are liftable. This proves theorem 2.14.

Open problem: give a direct proof of theorem 2.14 by writing down an
explicit lifting of $[a_1,b_1]$ as a non-commutative polynomial
in the $a_i$'s and the $b_i$'s. This seems surprisingly
difficult, possibly because such a lifting is far from unique so
it is hard to think of a canonical way to define it.

\proclaim Theorem 2.15. Suppose that $U$ is a torsion free bialgebra
over a principal ideal domain $R$, and suppose that the Lie algebra of
primitive elements of $U$ is a free $R$-module. If $U$ is
generated as an algebra
by all coefficients $a_i$ of some set of group-like elements
$1+a_1x+a_2x^2+\cdots\in U[[x]]$ then $U$ is a Hopf algebra with a
structural basis.

Proof.
As $R$ is a principal ideal domain, any submodule of a free
module is free, and in particular the Lie algebra $L$ of liftable
primitive elements of $U$ is a free $R$ module.
Suppose that $I$ is a totally ordered set indexing
a basis $a^i_1,i\in I$ of $L$. (The superscript on $a^i_1$ is
just used as an index and has nothing to do with powers of $a_1$.) For each
$a^i_1$ chose a lifting $a^i(x)=\sum_n a^i_n(x)$. Define
$H$ to be the $R$-module spanned by the elements
of the form $a^{i_1}_{n_1}\cdots a^{i_k}_{n_k}$ with $i_1<i_2\cdots <i_k$.

We need to check that the natural map from
$H\otimes_R H$ to $U\otimes_RU$ is injective,
so that $H\otimes_R H$ can be regarded as a submodule of $U\otimes_R U$.
As $U$ is torsion-free, its submodule $H$ is also torsion free
and therefore flat (as $R$ is a principal ideal domain).
Hence $H\otimes_R H\subseteq U\otimes_R U$ (as
$H\otimes_R H\subseteq U\otimes_RH\subseteq U\otimes_R U$).

The spanning set of $H$ of elements of the form
$a^{i_1}_{n_1}\cdots a^{i_k}_{n_k}$  satisfies the relations
of a structural basis, so in particular $H$ is a sub coalgebra of $U$.
This set of elements is linearly independent by lemma 2.7
because the primitive elements $a^i_1$ in it are linearly independent,
so it forms a structural basis for the coalgebra $H$.
We will show that $U$ has a structural basis by showing
that it is equal to $H$. To do this we will first
show that $H$ is closed under multiplication,
and then show that $H$ contains
a set of generators of the algebra $U$.

We define a filtration of $H$ by letting $H_n$ be the subspace
of $H$ spanned by all the monomials $a^{i_1}_{n_1}\cdots a^{i_k}_{n_k}$
with $n_1+\cdots n_k\le n$.
We will show by induction on $N$ that
$H_mH_n\subseteq H_{m+n}$ whenever $m+n\le N$; this is trivial for $N=1$.
To do this it is sufficient to show that for
any $r,s$ with $r+s\le N$ we have
$a^i_ra^j_s\equiv a^j_sa^i_r\bmod H_{r+s-1}$ whenever $i>j$, and
$a^i_ra^i_s\equiv {r+s\choose r} a^i_{r+s}\bmod H_{r+s-1}$.

We will show by induction on $k$ that if $1\le k\le N$ then
there is a group-like element $c(x,y)=\sum c_{ij}x^iy^j$ such that
$c_{ij}\in H_{i+j}$ if $i+j\le N$ and $c_{ij}= b_ja_i$ for
$i+j\le k$. This is obvious for $k=1$, as we can take $c_{ij}=a_ib_j$.
Suppose we have proved this for some value of $k<N$, so we
wish to prove it for $k+1$.
Consider the group-like element $b(y)a(x)c(x,y)^{-1}$. The coefficients
of $x^iy^j$ of this are 0 for $0<i+j\le k$, so by corollary
2.13 the coefficients for $i+j=k+1$ are all liftable. For each of
them choose a lifting $d_{i}(z)$, such that
the coefficient of $x^j$ of any of these liftings is in $H_j$.
We can construct such liftings
as products $a^{i_1}(r_1x)a^{i_2}(r_2x)\cdots$
if $d_{i} = r_1a^{i_1}+r_2a^{i_2}+\cdots$ with $i_1<i_2<\cdots$.
Then the product
$$e(x,y)=\Big(\prod_{1\le i\le k}d_{i}(x^iy^{k+1-i})\Big) c(x,y)$$
satisfies the conditions above for $k+1$. In fact by construction
we see that $e_{ij}=b_ja_i$ for $i+j\le k+1$. Also the coefficient
of $x^iy^j$ in $e(x,y)$ is either already an element of
$H_{i+j}$ or it is a sum of products of elements in $H_{i_1}, H_{i_2},\ldots$
with $i_1+i_2+\cdots <i+j$. If $i+j\le N$ then by the inductive assumption
on $N$ such  product is already in $H_{i+j-1}\subseteq H_{i+j}$, 
so the coefficient
of $x^iy^j$ is in $H_{i+j}$ provided $i+j\le N$.
This proves the inductive hypothesis for all $k\le N$.

If we take $k=N$ in the inductive hypothesis above, we now see that
$b_ja_i\in H_{i+j}$ for $i+j=N$. By induction on $N$
we see that $b_ja_i\in H_{i+j}$ for all $i,j$.

The proof that $a_ia_j={i+j\choose j}a_{i+j}\bmod H_{i+j-1}$
is similar, except that we use the group-like element
$a(x+y)=\sum_{i,j}{i+j\choose i}a_{i+j} x^iy^j$ instead
of $b(x)a(y)$.

This shows that $H$ is a subalgebra of $U$ (and hence a Hopf algebra with
a structural basis).

Finally we check that $H$ contains the coefficients
of any lifting $a(x)$ of any primitive element $a_1$ of $U$.
We will show by induction on $n$ that if $n\ge 1$ then
there is a lifting $b(x)$ of $a_1$ such that all coefficients of
$b(x)$ are in $H$ and $b_i=a_i$ for $i\le n$. This is clear for
$n=1$ as we just take any lifting of $a_1$ with coefficients in $H$.
Suppose that $b_i=a_i$ for $i\le n$. Then
$b(x)^{-1}a(x)$ has all its coefficients of $x^i$ vanishing
for $0<i\le n$, so by corollary 2.13 (with $k=1$) we see that
the coefficient of $x^{n+1}$ is liftable. Choose a lifting
$c(x)$ with coefficients in $H$. Then $b(x)c(x^{n+1})$ has all coefficients
in $H$ and the coefficients of $x^i$ for $i\le n+1$
are equal to $a_i$. This proves the inductive hypothesis for
$n+1$, and hence for all $n$. Therefore all coefficients
of $a(x)$ are in $H$.

We have shown that $H$ is a bialgebra with a structural basis
and that $H$ contains a set of generators of the algebra $U$.
Therefore $H$ is equal to $U$, and $U$ has a structural basis.
This proves theorem 2.15.

{\bf Example 2.16.} Suppose $n$ is a
nonnegative integer. The $\Z$-Hopf algebra $F_n$ is defined to the free
associative algebra over $\Z$ 
on a set of elements $Z_\alpha$, $\alpha\in \Z_{\ge 0}^n$,
$\alpha\ne 0$, with the comultiplication defined by
$$\Delta(Z_\alpha) =
\sum_{0\le \beta\le \alpha}Z_\beta\otimes Z_{\alpha-\beta}$$
on the generators (where $Z_0=1$).
Then it follows immediately from theorem 2.15 that $F_n$ has
a structural basis.
The fact that $F_n\otimes \F_p$ has a structural basis was first proved by
Dieudonn\'e; more precisely it follows from the comments at the end of
section 15 of [D], which show that $F_n(\F_p)$ has a structural basis,
together with theorem 3 of [D], which shows that
Dieudonn\'e's definition of $F_n(\F_p)$ is equivalent to $F_n\otimes
\F_p$. 
There are proofs that $F_1$ has a structural basis in [D72,
proposition 2.2] and [Sh], but unfortunately the paper [Sh] seems to
be unpublished, and the proof in [D72] has a gap (see [H p. 516]): in
[D72, page 5, line 6], it is implicitly assumed that $I_k/I_k^2$ is
torsion free. (It follows from example 2.16 that this is indeed true,
but it seems rather hard to prove.)

{\bf Example 2.17.}
Suppose that $U$ is a torsion free bialgebra
over a principal ideal domain $R$, and suppose that the Lie algebra of
primitive elements of $U$ is a free $R$-module. Then $U$ has
a subalgebra $H$ with a structural basis containing all
subalgebras with a structural base. The primitive elements of
$H$ are exactly the liftable primitive elements of $U$, and $H$ is generated
as an algebra by all coefficients of all liftings of primitive
elements of $U$. This follows
easily from theorem 2.15.

\proclaim 3.~Liftings of Lie algebra elements.

In this section we construct liftings of certain elements of the Lie
algebra of the vertex algebra of a double cover of an even lattice
([B86]).  More precisely we construct liftings of elements in the root
spaces of vectors of norms 2 or 0. This will be enough for
applications to the fake monster Lie algebra, because its rational
form is generated by the root spaces of roots of norms 2 or 0. For
vectors in norm 2 root spaces we construct a formal one parameter
group lifting any element in the root space. For norm 0 roots we
cannot do this, but have to make do with lifting vectors to formal
curves.

\proclaim Lemma 3.1. Suppose $a$ is an element of a vertex algebra,
and define $a^k$ for $k\ge 0$ by $a^0=1$, $a^{k+1} = a_{-1}(a^k)$.
Then
$$a^k_{n}(b) =
\sum_{0\le j\le k}{k\choose j}
\sum_{{i_1,\ldots,i_j<0\atop i_{j+1},\ldots, i_k\ge 0}\atop i_1+\cdots+i_k=n-k}
a_{i_1}\cdots a_{i_k}b
$$

Proof. This follows by induction on $n$ using the fact that
$$(a_{-1}a^k)_nb = \sum_{i\ge 0} a_{-1-i}(a^k_{n+i}b)
+\sum_{i\ge 0} a^k_{n-1-i}(a_ib)
$$
which is a special case of the vertex algebra identity.
This proves lemma 3.1.

\proclaim Corollary 3.2. Suppose that $a$ is an element
of a vertex algebra $V$
such that $a_ia=0$ for $i\ge -1$.
Then for each $k\ge 0$ the operator
$a_0^{k}$ is divisible by $k!$.

Proof. We show 
that $a_0^k$ is divisible
by $k!$ by induction on $k$, which is obvious for $k=0,1$. 
As $a_{-1}a=0$ the left hand side of
lemma 3.1 vanishes for $k\ge 2$.
As $a_ia=0$ for $i\ge 0$ all the operators $a_i$
commute with each other
because of the formula
$$[a_i,b_j] = \sum_{k\ge 0} {i\choose k} (a_kb)_{i+j-k},$$
which together with the induction
hypothesis implies that the sum of all terms on the
right hand side of lemma 3.1 other than $a_0^k$ is divisible by
$k!$. Hence $a_0^k$ is also divisible by $k!$. This proves lemma 3.2.

\proclaim Corollary 3.3. Suppose $V$ is the vertex algebra
(over $\Z$) of a double cover $\hat L$ of an even integral lattice $L$.
If $\alpha\in \hat L$ has positive norm and $a$ is the element
$e^\alpha$ of $V$, then all coefficients of
$$\exp(xa_0)=\sum_kx^ka_0^k/k!$$ map $V$ into $V$.

Proof. For any $\alpha, \beta\in L$ we have
$e^\alpha_ie^\beta = 0$ for $i+(\alpha,\beta)\ge 0$.
In particular $e^\alpha_i e^\alpha=0$
for $i\ge -1$ if $(\alpha,\alpha)>0$. 
Corollary 3.3 now follows from corollary 3.2.

Corollary 3.3 will allow us to lift elements in the root spaces of
norm 2 vectors. In the rest of this section we show how to 
lift elements in the root spaces of norm 0 vectors. 

\proclaim Lemma 3.4.
Suppose that the the elements $\Gamma_i$ for $i\in \Z$ are independent
formal variables. If $I=(i_1,\ldots, i_m)$ is a finite sequence of
integers then define $\Gamma_I$ to be $\Gamma_{i_1}\cdots
\Gamma_{i_m}$ and define $\Sigma(I)$ to be $i_1+\cdots +i_m$ and
define $l(I)$ to be $m$. Then all coefficient of the power series
$$E=\exp\Big( \sum_{m,n>0}
\sum_{I\in \Z^m}
\sum_{J\in \Z^n\atop \Sigma(J)=-\Sigma(I)>0}
{\Sigma(I)\over mn}\Gamma_{I}\Gamma_J\Big)$$
are integers.

Proof.  If $I$ is any finite sequence of integers then $I^k$ means the
obvious concatenation of $k$ copies of the sequence $I$.  We will call
a pair of finite sequences primitive if it is not of the form
$(I^k,J^k)$ for some $k\ge 2$.  Any pair of finite sequences can be
written uniquely as $(I^k,J^k)$ for some primitive pair $(I,J)$,
which we call the primitive core of the pair $(I^k,J^k)$.  Consider
the group $\Z\times \Z$ acting on pairs $(I,J)$ by the first $\Z$
acting as cyclic permutations of the elements of $I$ and the second
$\Z$ acting as cyclic permutations of the elements of $J$.  We group
the pairs $(I,J)$ indexing the terms of the sum in the exponent of $E$
into equivalence classes, where we say two pairs of sequences are
equivalent if their primitive cores are conjugate under $\Z\times \Z$.
For any primitive element $(I,J)$ with $\Sigma(I)+\Sigma(J)=0$ we will
show that all coefficients of the exponential of the sum of all terms
in the equivalence class of $(I,J)$ are integral. This will show that
the coefficients of $E$ are integral because $E$ is an infinite
product over the set of orbits of primitive elements of expressions
like this.  Let $m'$ and $n'$ be the number of orbits of $I$ and $J$
under the cyclic action of $\Z$.  Then the exponential of the sum of
the terms equivalent to $(I,J)$ is
$$\eqalign{
&\exp\Big(\sum_{k>0}{\Sigma(I^k)\over l(I^k)l(J^k)}
(\hbox{number of orbits of $(I,J)$ under $\Z\times\Z$})
 \Gamma_{I^k}\Gamma_{J^k}\Big)\cr
=&\exp\Big(\sum_{k>0} {k\Sigma(I)m'n'\over kmkn}
\Gamma_{I^k}\Gamma_{J^k}\Big)\cr
=&\exp\Big({\Sigma(I)m'n'\over mn} \sum_{k>0} \Gamma_I^k\Gamma_J ^k/k\Big)\cr
=&\exp\Big(-{\Sigma(I)m'n'\over mn} \log(1- \Gamma_I\Gamma_J )\Big)\cr
=& (1- \Gamma_I\Gamma_J )^{-\Sigma(I)m'n'/mn}\cr
}$$
The number $\Sigma(I)$ is divisible by $m/m'$ as $I$ is the concatenation
of $m/m'$ identical sequences and $\Sigma(I)$
is the sum of the elements of $I$,
and similarly $\Sigma(I)=-\Sigma(J)$ is divisible
by $n/n'$. Moreover $m/m'$ and $n/n'$ are coprime as $(I,J)$ is primitive.
Hence $\Sigma(I)m'n'/mn$ is an integer because it is equal to $\Sigma(I)$
divided
by two coprime factors $m/m'$ and $n/n'$ of $\Sigma(I)$.
This implies that $(1- \Gamma_I\Gamma_J )^{-\Sigma(I)m'n'/mn}$ has
integral coefficients,
and so $E$ does as well as it is an infinite product of expressions like this.
This proves lemma 3.4.

\proclaim Theorem 3.5. Let $V$ be the vertex algebra
(over the integers) of the double cover of an even lattice $L$.
Define $D^*V$ to be the sum of the spaces $D^{(i)}V$ for $i\ge 1$.
We recall from [B86] that $V/D^*V$ has a natural Lie algebra structure,
with the bracket defined by $[u,v]=u_0v$. Moreover this Lie algebra
acts on $V$ preserving the vertex algebra structure of $V$.
Suppose that $\alpha, \gamma\in L$ with
$\alpha$ orthogonal to the norm 0 vector $\gamma$.
Then
$$\exp\Big(\sum_{i>0} {x^i\over i}(\alpha(1)e^{i\gamma})_0 \Big)$$
is a lifting of $(\alpha(1)e^\gamma)_0$ in the universal enveloping algebra
of $(V/D^*V)\otimes \Q$ all of whose
coefficients map $V$ to $V$.

Proof. It is obvious that the element is group-like as it is
the exponential of a primitive element, so the only problem
is to show that it preserves the integral form $V$.

This element is also an automorphism of the vertex algebra $(V\otimes
\Q)[[x]]$ as it is the exponential of a derivation of this vertex
algebra. Hence to show its coefficients
preserve the integral form of $V$ 
it is sufficient to show that it maps each of the generators
$e^\beta$, $\beta\in L$, of the vertex algebra $V$ into $V[[x]]$.

Define the operators $\Gamma_i$ by
$$e^\gamma(z)= \sum_i\Gamma_iz^i.$$
Recall the following formulas from [B86].
$$\alpha(1)(z) = \sum_j \alpha(1)_{-j}z^{j-1} = \sum_j \alpha(j)z^{j-1}$$
All coefficients of $e^\gamma$ commute with everything in sight
as $\gamma$ is orthogonal to $\gamma$ and $\alpha$.
So $(\alpha(1)e^{i\gamma})_0$ is equal to
the coefficient of $z^{-1}$ in
$$\alpha(1)(z)(e^\gamma(z))^i
=\sum_j \alpha(j)z^{j-1} \sum_{I\in \Z^i}\Gamma_Iz^{\Sigma(I)}$$
Therefore $\sum_{i>0}{x^i\over i}(\alpha(1)e^{i\gamma})_0 $ is equal to
$A^++A^0+A^-$ where
$$\eqalign{
A^+&= \sum_{m>0} \sum_{j>0} \alpha(j) \sum_{I\in \Z^m\atop \Sigma(I)=-j}
{x^m\over m} \Gamma_I \cr
A^0&= \sum_{m>0} \sum_{j=0} \alpha(j) \sum_{I\in \Z^m\atop \Sigma(I)=-j}
{x^m\over m} \Gamma_I \cr
A^-&= \sum_{m>0} \sum_{j<0} \alpha(j) \sum_{I\in \Z^m\atop \Sigma(I)=-j}
{x^m\over m} \Gamma_I \cr
}$$

We would like to pull out a factor of $\exp(A^-)$
from $\exp(A^++A^0+A^-)$ because $\exp(A^-)(e^\beta)
=e^\beta$, but we have to be careful when doing this because $\alpha(i)$
does not commute with $\alpha(-i)$ if $i\ne 0$.

We can easily evaluate $[A^+,A^-]/2$ using the fact that
$[\alpha(i),\alpha(j)]=j(\alpha,\alpha)$ if $i+j=0$ and 0 otherwise, and
we find that
$$[A^+,A^-]/2= {1\over 2} (\alpha,\alpha)
\sum_{m,n>0}\sum_{k<0}
\sum_{I\in \Z^m\atop \Sigma(I)=k}
\sum_{J\in \Z^n\atop \Sigma(J)=-k}
{k\over mn}\Gamma_I\Gamma_J.$$
In particular this verifies the fact used below that $[A^+,A^-]$ commutes with
$A^+$ and $A^-$.
The exponential of this is the expression in lemma 3.4 raised to the power
of $(\alpha,\alpha)/2$, and therefore has integral coefficients by
lemma 3.4 and because $(\alpha,\alpha)$ has even norm.
So $\exp([A^+,A^-]/2)$ maps $V$ into $V[[x]]$.

Recall the
formula
$$\exp(A^++A^-)=\exp([A^-,A^+]/2)\exp(A^+)\exp(A^-)
$$
which is valid because $[A^+,A^-]$ commutes with $A^+$ and $A^-$.
(This is essentially just the first few terms of the
Baker-Campbell-Hausdorff formula.)  Now we look at
$$\exp(A^++A^0+A^-)e^\beta =
\exp([A^-,A^+]/2)\exp(A^+)\exp(A^0)\exp(A^-)e^\beta
$$
It is obvious that $\exp(A^-)(e^\beta) = e^\beta$,
and we have checked above that all coefficients of
$\exp([A^-,A^+]/2)$ map $V$ to $V$. Hence to complete the
proof of theorem 3.5 it is sufficient to prove that
all coefficients of $\exp(A^0)$ and $\exp(A^+)$ map $V$ into $V$.

We check that $\exp(A^0)$ maps $V$ to $V[[x]]$. The follows
because $A^0$ is an infinite sum of expressions
like
$$\sum_{n>0} \alpha(0) \Gamma_{J^n}x^n/n=-\alpha(0)\log(1-x\Gamma_J)$$
where $J$ is a primitive sequence (in other words a sequence
that cannot be written in the form $I^m$ for some $m>1$).
Hence $\exp(A^0)$ is an infinite product of
terms of the form $(1-x\Gamma_{J})^{-\alpha(0)}$, which
map $V[[x]]$ to itself because $\alpha(0)$ has integral eigenvalue
$(\alpha,\beta)$ on the subspace of $V$ of degree $\beta\in L$.
This shows that all coefficients of $\exp(A^0)$ map $V$ into $V$.

Finally we have to show that all coefficients of $\exp(A^+)$ map $V$
into $V$.
As usual we divide the sum over elements $I$ in $A^+$ into
classes consisting of powers of conjugates of primitive elements $I$.
We see that $A^+$ is a sum over all orbits of primitive elements $I$
with $\Sigma(I)>0$ of
expressions like
$$l(I)\sum_{k>0}\alpha(\Sigma(I^k)){x^{l(I^k)}\over l(I^k)}\Gamma_{I^k}
=\sum_{k>0}\alpha(k\Sigma(I)){x^{kl(I)}\over k} \Gamma_I^k.
$$
(The factor at the front is the number
of conjugates of $I$ under the cyclic action of $\Z$, which
is equal to $l(I)$ because $I$ is primitive.)
So it is sufficient to show that the exponential of this
expression has integral coefficients.
Let $y_1,y_2,\ldots, $ be a countable number of independent
variables,
and identify $\alpha({k})$ with the symmetric function $\sum y_i^k$
of the $y$'s for $k>0$. (See [M chapter 1].)
Then
$$\eqalign{
&\exp\Big(\sum_{k>0}\alpha(k\Sigma(I)){x^{kl(I)}\over k} 
\Gamma_I^k\Big)\cr
=&\exp\Big(\sum_i\sum_{k>0}y_i^{k\Sigma(I)}{x^{kl(I)}\over k} 
\Gamma_I^k\Big)\cr
=&\exp\Big(\sum_i -\log(1-y_i^{\Sigma(I)}{x^{l(I)}} 
\Gamma_I)\Big)\cr
=&\prod_{ i}{1\over 1-y_i^{\Sigma(I)}{x^{l(I)}} 
\Gamma_I }.\cr
}$$
So we see that
$$\exp(A^+) = \prod_I\prod_{ i}{1\over 1-y_i^{\Sigma(I)}{x^{l(I)}}
\Gamma_I }
$$ 
where the product over $I$ is a product over all orbits
of primitive sequences $I$ with $\Sigma(I)>0$ under the cyclic action
of $\Z$.  The last line is a power series in the elements $\Gamma_I$
and $x$ whose coefficients are symmetric functions in the $y$'s, and
hence are polynomials with integral coefficients in the complete
symmetric functions of the $y$'s. So we have to show that each
complete symmetric function, considered as a polynomial in the
$\alpha$'s with rational (not necessarily integral!) coefficients maps
$V$ to $V$.  The complete symmetric functions are the polynomials
$e^{-\alpha}D^{(n)}(e^\alpha)$ considered as elements of the ring $V$
underlying the vertex algebra $V$ (as follows from [M, Chapter 1,
2.10]). By definition of the integral form $V$ these polynomials map
$V$ to itself.

This proves theorem 3.5.

\proclaim 4.~The fake monster smooth Hopf algebra.

We recall the construction of the fake monster Lie algebra [B90].
It is the Lie algebra of physical states of the vertex algebra of
the double cover of the lattice $II_{25,1}$. This Lie algebra has
an integral form $\m$ consisting of the elements represented by
elements of the integral form of the vertex algebra. 
We recall some properties of $\m$:
\item{1} $\m$ is graded by the lattice $II_{25,1}$,
and the piece $\m_\alpha$ of degree $\alpha\in II_{25,1}$ has dimension
$p_{24}(1-\alpha^2/2)$ if $\alpha\ne 0$ and 26 if $\alpha=0$,
where $p_{24}(n) $ is the number of partitions of $n$ into
parts of 24 colors.
\item{2} $\m$ has an involution $\omega$ lifting the
involution $-1$ of $II_{25,1}$.
\item{3} $\m$ has a symmetric invariant integer valued bilinear
form $(,)$, and the pairing between $\m_\alpha\otimes \Q$
and $\m_{-\alpha}\otimes \Q$ is nonsingular.
\item{4} $\m\otimes \Q$ is a generalized Kac-Moody algebra.
The simple roots are given by the norm 2 vectors $r$
with $(r,\rho)=-1$, together with all positive multiples of $\rho$
with multiplicity 24, where $\rho$ is a primitive norm 0 vector
of $II_{25,1}$ such that $\rho^\perp/\rho$ is isomorphic to the Leech lattice.
(This follows from [B90. theorem 1].)

We define $U^+(\m)$ to be the $\Z$-subalgebra of
the universal enveloping algebra $U(\m\otimes \Q)$ generated
by the coefficients of the liftings of elements in root spaces of
the simple roots and their negatives
constructed in corollary 3.3 and theorem 3.5. (Note that all simple roots
of $\m$
have norms 2 or 0 so we can always apply one of these two types of liftings.)

\proclaim Theorem 4.1. There is a $II_{25,1}$-graded
Hopf algebra $U^+(\m)$ over $\Z$ with the following
properties.
\item {1} $U^+(\m)$ has a structural basis over $\Z$.
\item {2} The primitive elements of $U^+(\m)$ are an integral form of the
fake monster Lie algebra $\m$.
\item {3} For every norm 2 vector of $II_{25,1}$, $U^+(\m)$ contains
the usual (Kostant) integral form of the universal enveloping algebra
of the corresponding $sl_2(\Z)$.

Proof.  The algebra $U^+(\m)$ is a $\Z$-Hopf subalgebra of
$U(\m\otimes\Q)$ as it is generated by coalgebras. Also $U^+(\m)$ is
obviously torsion free as it is contained in a rational vector space.
It is easy to check directly that the degree zero primitive elements
of $U^+(\m)$ are just the degree 0 elements of $\m$ and therefore form
a free $\Z$ module (of rank 26).  If $\alpha$ is any nonzero element
of $II_{25,1}$ having nonzero inner product with some element
$\beta\in II_{25,1}$, then  $(\alpha,\beta)u\in
\m_\alpha$ for any primitive element $u\in U^+(\m)$ of degree
$\alpha$, because $U^+(\m)$ maps $\m$, and hence $\beta$, to
$\m$ by corollary 3.3 and theorem 3.5. 
This shows that all root spaces of primitive elements of
$U^+(\m)$ are free $\Z$-modules.  We can now apply theorem 2.15 to see
that $U^+(\m)$ is a Hopf algebra with a structural basis.  The
primitive elements of $U^+(\m)$ form an integral form of the fake
monster Lie algebra, because the fake monster Lie algebra over the
rationals is generated by the root spaces of simple roots and their
negatives, and all simple roots have norms 2 or 0.  This proves theorem 4.1.

\proclaim 5.~A smooth Hopf algebra for the Virasoro algebra.

In this section we show that there is a Hopf algebra over $\Z$ with a
structural basis whose primitive elements form the natural
integral form of the
Virasoro algebra (theorem 5.7). In other words, there is a formal
group law over $\Z$ corresponding to the Virasoro algebra. Moreover this
Hopf algebra acts on the integral form of the vertex algebra of any
even self dual lattice.

Let $R$ be a commutative ring.  We write $Hom(R,R)$ for the ring of
homomorphisms of the abelian group $R$ to itself, and $Der(R)$ for the
Lie algebra of derivations of the ring $R$, and $U(Der(R))$ for the universal
enveloping algebra of $Der(R)$ over $\Z$.  Consider the group of all
element $a=\sum_i a_i\epsilon^i\in Hom(R,R)[[\epsilon]]$ with $a_0=1$
that induce automorphisms of the $\Z[[\epsilon]]$ algebra
$R[[\epsilon]]$.  We can think of the elements of this group
informally as ``infinitesimal curves in the group of automorphisms of
$Spec(R)$''.  We will call a derivation of $R$ liftable if it is of
the form $a_1$ for some $a$ as above.

\proclaim Lemma 5.1. Let  $R$ be any commutative algebra
with no $\Z$-torsion. Then any $a\in  Hom(R,R)[[\epsilon]]$
with $a(0)=1$ that is an automorphism  of  $R[[\epsilon]]$
is the image of a unique group-like element $G_a$ of
$U(Der(R)\otimes \Q)[[\epsilon]]$ with $G_a(0)=1$ under
the natural map from $U(Der(R)\otimes \Q)$ to $Hom(R,R)\otimes \Q$.

Proof. We note that $\log(a)$
is a well defined element of $(Hom(R,R)\otimes \Q)[[\epsilon]]$
as $a=1+O(\epsilon )\in Hom(R,R)[[\epsilon ]]$.
As $a$ is a ring homomorphism it follows that
that $\log(a)$ is a derivation of $(R\otimes \Q)[[\epsilon]]$
and is therefore an element of $(Der(R)\otimes \Q)[[\epsilon]]$.
Now we consider $\log(a)$
to be an element of the universal enveloping algebra
$U(Der(R)\otimes \Q)[[\epsilon ]]$
and we define $G_a\in U(Der(R)\otimes \Q)[[\epsilon ]]$ by
$$G_a = \exp(\log(a)),$$
where the log is computed in $(Hom(R,R)\otimes \Q)[[\epsilon]]$
and the exponential
is computed in $U(Der(R)\otimes \Q)[[\epsilon ]]$.
It is obvious that the action of $G_a$ on $R[[\epsilon]]$ is the same
as that of $a$.
Also $G_a$ is group-like because it is the exponential of a primitive
element. It is easy to check that $G_a$ is the unique group-like
lifting of $a$ with $G_a(0)=1$, because the log of a group-like element
must be primitive and must therefore be the same as $\log(a)$.
This proves lemma 5.1.

\proclaim Corollary 5.2. Suppose that $R$ is a commutative
ring with no $\Z$ torsion such that $Der(R)$ is a free
$\Z$-module. Define $U^+(Der(R))$ to be the subalgebra of
$U(Der(R\otimes\Q))$ generated by all the coefficients of all
group-like elements of $U(Der(R)\otimes \Q)[[\epsilon]]$ that have
constant coefficient 1 and map $R[[\epsilon]]$ to $R[[\epsilon]]$. Then
$U^+(Der(R))$ is a Hopf algebra over $\Z$ with a structural basis, and
its primitive elements are the liftable primitive elements of
$Der(R)$.

Proof. Applying theorem 2.15 shows that $U^+(Der(R))$
is a Hopf algebra over $\Z$ with
a structural basis. By lemma 5.1 the space of primitive elements
of $U^+(Der(R))$
is the same as the space of liftable primitive elements of
$Der(R)$. This proves corollary 5.2.

{\bf Example 5.3} Suppose we take $R$ to be the algebra
$\Z[x][x^{-1}]$ of Laurent polynomials.  Then $Der(R)=Witt$ is the
Witt algebra over $\Z$, which is spanned by the elements $L_m=-x^{m+1}{d\over
dx}$ for $m\in \Z$.  All elements $L_m$ are liftable; for example, we
can use the automorphism of $R[[\epsilon]]$ taking $x$ to $x-\epsilon
x^{m+1}$ to show that $L_m$ is liftable. Therefore the Hopf algebra
$U^+(Witt)$ is a Hopf algebra over $\Z$ with a structural basis, whose
primitive elements are exactly the elements of the Witt algebra.

\proclaim Lemma 5.4. Let $N$ be an integer. Let $R_N$ be the representation
of $Witt$ with a basis of elements $e_n$, $n\in \Z$, with the action
given by $L_m(e_n)= (Nm+n)e_{m+n}$. Then the action of
$Witt$ on $R_N$ can be extended to an action of $U^+(Witt)$ on $R_N$.

Proof. For $N=-1$ the module $R_N$ is the module of first order
differential operators on $R=\Z[x][x^{-1}]$, so the automorphisms $G_a$
extend to $R_N[[\epsilon]]$. For other negative values of $N$ the $R$-module 
$R_N$
is a tensor product of $-N$ copies of the module $R_{-1}$, and
$R_N$ for $N$ positive is the dual of $R_N$ for $N$ negative.
Therefore $U^+(Witt)$ extends to these modules as well. This proves
lemma 5.4.

We define $Witt_{\ge n}$ for $n=-1,1$ to be the subalgebra of
$Witt$ spanned by $L_i$ for $i\ge n$.
Let $V$ be the vertex algebra of some even lattice. It contains elements
$e^\alpha$ for $\alpha\in L$, so there are operators
$e^\alpha_i$ on $V$ for $i\in \Z$. The algebra $Witt_{\ge -1}$ also
acts naturally on $V$. 

\proclaim Lemma 5.5. 
Put $N=\alpha^2/2-1$ and $e_j=e^\alpha_{\alpha^2/2-1-j}$ for $j\in \Z$
and let $R_N$ 
be the space with the elements $e_j$ as a basis. Define an action of 
the algebra $U^+(Witt_{\ge -1})$ on $R_N\otimes U^+(Witt_{\ge -1})$
using the action  on $R_N$ as in lemma 5.4 and the action on
$U^+(Witt_{\ge -1})$ by left multiplication and the coalgebra structure
of $U^+(Witt_{\ge -1})$.  If $u\in U^+(Witt_{\ge -1})$ then
$ue^\alpha_m= \sum_i e^\alpha_i u_i$ as operators on $V$, where $\sum_i
e^\alpha_i\otimes u_i=u(e^\alpha_m\otimes 1)$ is the image of
$e^\alpha_m\otimes 1$ under the action of $u$ on $R_N\otimes
U^+(Witt_{\ge -1})$. In particular if $U^+(Witt_{\ge -1})$
maps some element $v\in V$ to $V$ then it also maps $e^\alpha_i(v)$ to $V$
for any $\alpha\in L$ and $i\in \Z$.

Proof.
If $u$ is of the form $L_i$ for $i\ge -1$ this can be proved as follows.
A standard vertex algebra calculation shows that
$$[L_i, e^\alpha_j] = ((i+1)(\alpha^2/2-1)-j) e^\alpha_{i+j}.$$ 
This
shows that lemma 5.5 is true when $u\in Witt_{\ge -1}$.  If lemma 5.5
is true for two elements $u$, $u'$ of $U^+(Witt_{\ge -1})$ then it is
true for their product. If it is true for some nonzero integral
multiple of $u\in U^+(Witt_{\ge -1})$ then it is true for $u$ because
$V$ is torsion-free.  To finish the proof we observe that the algebra
$U^+(Witt_{\ge -1})$ is generated up to torsion by the elements $L_m$
for $m\ge -1$.  This proves lemma 5.5.

\proclaim Lemma 5.6. Suppose that $V$ is the vertex algebra
of the double cover of some even lattice, with the standard action of
$Witt_{\ge-1}$ on $V\otimes \Q$ ([B86]).  This action extends to an action of
$U^+(Witt_{\ge -1})$ on $V$.

Proof. The vertex algebra $V$ is generated from the element $1$ by the
actions of the operators $e^\alpha_n$ for $n\in \Z$ and $\alpha\in L$.
Define $F^n(V)$ by defining $F^0(V)$ to be the space spanned by 1, and
defining $F^{n+1}(V)$ to be the space spanned by the actions of
operators of the form $e^\alpha_m$ on $F^n(V)$. Then $V$ is the union
of the spaces $F^n(V)$, so it is sufficient to prove that each space
$F^n(V)$ is preserved by the action of $U^+(Witt_{\ge -1})$. We will
prove this by induction on $n$. For $n=0$ is is trivial because
$L_n(1)=0$ for $n\ge -1$. If it is true for $n$, then it follows
immediately from lemma 5.5 that it is true for $n+1$.
This proves lemma 5.6.

The algebra $U^+(Witt_{\ge -1})$ can be $\Z$-graded in such a way that
$L_m$ has degree $m$. This follows easily from the fact that
we can find graded liftings of the elements $L_m$. We define
$U^+(Witt_{\ge 1})$ to be the subalgebra generated by the coefficients
of graded liftings of the elements $L_m$ for $m\ge 1$. It is easy to check
that this is a Hopf algebra with a structural basis whose Lie algebra of
primitive vectors has a basis consisting of the elements $L_m$ for $m\ge 1$.
It is $\Z_{\ge 0}$-graded, with all graded pieces being finite dimensional;
in fact the piece of degree $n\in \Z$ has dimension $p(n)$ where
$p$ is the partition function.

We let $U(Witt_{\ge 1})$ be the
universal enveloping algebra of the Lie algebra $Witt_{\ge 1}$.
It is $\Z_{\ge 0}$-graded in the obvious way and is a subalgebra
of $U^+(Witt_{\ge 1})$.

In the rest of this section we construct an integral form with a
structural basis for the universal enveloping algebra of the Virasoro
algebra over $\Z$. This result is not used elsewhere in this paper.
We recall that the Virasoro algebra $Vir$ is a central extension
of $Witt$ and is spanned by elements $L_i$ for $i\in \Z$ and
an element $c/2$ in the center, with
$$[L_m,L_n]=(m-n)L_{m+n} +{m+1\choose 3}{c\over 2}.$$
We identify $Witt_{\ge -1}$ with the subalgebra of $Vir$
spanned by $L_m$ for $m\ge -1$.
The Virasoro algebra $Vir$ has an automorphism
$\omega$ of order 2 defined by $\omega(L_m)=-L_{-m}$, $\omega(c)=-c$,
and $\omega$ extends to an automorphism of the universal enveloping algebra
$U(Vir\otimes\Q)$. We define $U^+(Vir)$ to be the
subalgebra of $U(Vir\otimes\Q)$ generated by
$U^+(Witt_{\ge -1})$ and $\omega(U^+(Witt_{\ge -1}))$.

For any even integral lattice $L$ there is a double cover
$\hat L$, unique up to non-unique isomorphism,
such that $e^ae^b=(-1)^{(a,b)}e^be^a$.
We let $V_{\hat L}$ be the (integral form of the)
vertex algebra of $\hat L$. This is $L$-graded, and
has a self dual bilinear form on it (more precisely,
each piece of given $L$-degree $\alpha$ and given eigenvalue under $L_0$
is finite dimensional and dual to the piece of degree $-\alpha$),
and if $L$ is self dual then
$V_{\hat L}$
has a conformal vector generating an action of the Virasoro algebra.

\proclaim Theorem 5.7. The subalgebra $U^+(Vir)$ of
$U(Vir\otimes \Q)$ has the following properties:
\item 1 $U^+(Vir)$ is a $\Z$-Hopf algebra with a structural basis.
\item 2 The Lie algebra of
primitive elements of $U^+(Vir)$ has a basis consisting of
the elements $L_n$ for $n\in \Z$ and the element $c/2$.
\item 3 $U^+(Vir)$ maps the vertex algebra (over $\Z$) of any even
self dual lattice to itself.

Proof. We first construct the action on the vertex algebra $V$
of an even self dual lattice. We have an action of $U^+(Witt_{\ge -1})$
on $V$ by lemma 5.6. The vertex algebra of an even self dual lattice
is also self dual under its natural bilinear form, so the adjoint
of any linear operator on $V$ is also a linear operator on $V$.
The adjoint of $L_m$ is $L_{-m}$, so the adjoint $U^+(Witt_{\le 1})=
\omega(U^+(Witt_{\ge -1}))$ also maps $V$ to itself. As these two
algebras generate $U^+(Vir)$ this proves that $U^+(Vir)$ acts on $V$.

Next we find the Lie algebra $P$ of primitive elements of $U^+(Vir)$.
The elements $L_{m}$ for $m\ge -1$ are obviously in $U^+(Vir)$ because
they are in $U^+(Witt_{\ge -1})$, and similarly $L_m$ for $m\le 1$ is in
$P$. The element $c/2$ is in $P$ because $[L_2,L_{-2}]= 4L_0+c/2$.
So $P$ contains the basis described in theorem 5.7, and we have
to prove that $P$ contains no elements other than linear combinations of these.
As $U^+(Vir)$ and hence $P$ are both $\Z$-graded with $L_m$ having
degree $m$, it is sufficient to show that the degree $m$ piece of
$P$ is spanned by $L_m$ if $m\ne 0$ and by $L_0$ and $c/2$ if $m=0$.
For $m\ne 0$ this is easy to check as we just map the Virasoro algebra
to the Witt algebra and use example 5.3. The case $m=0$ is harder
and we will use the actions on vertex algebras of even self dual
lattices $L$ constructed above. If $L$ is such a lattice then
$c$ acts on $V$ as multiplication by $\dim(L)$, and $L_0$ has
eigenspaces with eigenvalue any given positive integer (at least
if $L$ has positive dimension). So if $xL_0+yc/2$ is in $P$ for some
$x,y\in \Q$ then $xm+yn/2$ is an integer whenever $m$ is a positive
integer and $n$ is the dimension of a nonzero even self dual lattice.
As we can find such lattices for any positive even integer $n$, this
implies that $x$ and $y$ are both integers. Hence
the degree 0 piece of $P$ is spanned by $L_0$ and $c/2$. This completes
the proof that $P$ is spanned by $L_m, m\in \Z$, and $c/2$.

Finally we have to show that $U^+(Vir)$ has a structural basis.
We know that all the elements $L_m$ for $m\ge -1$ are liftable in
$U^+(Witt_{\ge -1})$, so they are also liftable in $U^+(Vir)$.
Similarly the elements $L_m$ for $m\le 1$ are liftable in 
$\omega(U^+(Witt_{\ge -1}))$ and therefore in $U^+(Vir)$.
It is trivial to check that $c$ is liftable as it acts
as multiplication by some integer, so every primitive element of
$U^+(Vir)$ is liftable by lemma 2.5 and part 2 of theorem 5.7. 
Also $U^+(Vir)$ is generated
by the coefficients of group-like elements because this is true
for $U^+(Witt_{\ge -1})$ and its conjugate under $\omega$. The fact that
$U^+(Vir)$ has a structural basis now follows from theorem 2.15.
This completes the proof of theorem 5.7.

Example. If $L$ is an even lattice of odd dimension then $L_{-2}(1)$
is not in the vertex algebra of $L$, because $[L_2,L_{-2}] =
4L_0+\dim(L)/2$.  So part 3 of theorem 5.7 is false without the
assumption that $L$ is self dual.

\proclaim 6.~The no-ghost theorem over $\Z$.

The no-ghost theorem of Goddard and Thorn [G-T] implies that over
the reals, the contravariant form restricted to the
degree $\beta\ne 0$ piece of the fake monster Lie algebra $\m$
is positive definite and in particular nonsingular. We give a refinement of
this to the integral form of the degree $\beta$ piece $\m_\beta$, showing
that the any prime dividing the discriminant of the quadratic
form on $\m_\beta$ also divides the vector $\beta$.
In particular if $\beta$ is a primitive vector then 
$\m_\beta$ is a self dual positive definite integral lattice.
I do not know whether or not the discriminant can be divisible by $p$ when
$\beta$ is divisible by $p$.

If $\lambda=1^{i_1}2^{i_2}\cdots$ is a partition with $i_1$ 1's, $i_2$
2's, and so on, then we define $P(\lambda)$ to be the integer
$1^{i_1}2^{i_2}\cdots$ and $F(\lambda)$ to be $ i_1!i_2!\cdots$ and
$|\lambda|$ to be $1i_1+2i_2+\cdots$ and $l(\lambda)$ to be
$i_1+i_2+\cdots$. We also define $p(n)$ to be the number of partitions
of $n$.

\proclaim Lemma 6.1. Suppose that $n$ is an integer.
Then
$$\prod_{|\lambda|=n} P(\lambda)=\prod_{|\lambda|=n} F(\lambda).$$

Proof. We will show that both sides are equal to
$$\prod_{i>0} i^{\Sigma_{j>0}p(n-ij)}.$$

The left hand side $\prod_{|\lambda|=n} P(\lambda)$ is equal to
$\prod_{i>0} i^{n(i)}$ where $n(i)$ is the number of times $i$ occurs
in some partition of $n$, counted with multiplicities.  This number
$n(i)$ is equal to $\sum_j n(i,j)$ where $n(i,j)$ is the number of
times that $i$ occurs at least $j$ times in a partition of $n$. But
$n(i,j)$ is equal to $p(n-ij)$ because any partition of $n$ in which
$i$ occurs at least $j$ times can be obtained uniquely from a
partition of $n-ij$ by adding $j$ copies of $i$.  So the left hand
side $\prod_{|\lambda|=n} P(\lambda)$ is equal to $\prod_{i>0}
i^{\Sigma_{j>0}p(n-ij)}$.

On the other hand the right hand side $\prod_{|\lambda|=n} F(\lambda)$
is equal to $\prod_{i>0} i^{m(i)}$ where $m(i)$ is the number of times that
there is a partition of $n$ with some number occurring at least $i$ times
(counting a partition several times if it has more than
one number occurring at least $i$ times). But $m(i)$ is equal to
$\sum_{j}n(j,i)=\sum_{j}p(n-ji)$.
This shows that the right hand side
$\prod_{|\lambda|=n} F(\lambda)$ is also  equal to
$\prod_{i>0} i^{\Sigma_{j>0}p(n-ij)}$. This proves lemma 6.1.

\proclaim Lemma 6.2.
The submodule $U(Witt_{\ge 1})_n$ has index $\prod_{|\lambda|=n} F(\lambda)$
inside $U^+(Witt_{\ge 1})_n$.

Proof. Choose a graded lifting $1+a_{i,1}x+a_{i,2}x^2+\ldots$
of $L_i$ for each $i>0$. Then
the elements $a_{1,i_1}a_{2,i_2}\cdots$ for $1i_1+2i_2+\cdots=n$
form a base for $U^+(Witt_{\ge 1})$,
and the elements $i_1!a_{1,i_1}i_2!a_{2,i_2}\cdots$ for
$1i_1+2i_2+\cdots=n$ form a basis for $U(Witt_{\ge 1})_n$.
Therefore the index of $U(Witt_{\ge 1})_n$ in $U^+(Witt_{\ge 1})_n$
is
$$\prod_{i_1+2i_2+\cdots=n}i_1!i_2!\cdots=\prod_{|\lambda|=n}F(\lambda).$$
This proves lemma 6.2.

\proclaim Lemma 6.3.
Let $\gamma$ be a norm 0 vector of $II_{25,1}$. Suppose that $W$ is
the graded space of all elements generated by the action of the
elements $e^{-\gamma}D^{(i)}(e^\gamma)$ on $e^\beta$, so that $W$ is
acted on by the smooth integral form $U=U^+(Witt_{>0})$. Then the
graded dual $W[1/(\beta,\gamma)]^*$ of $W[1/(\beta,\gamma)]$ is a free
$U[1/(\beta,\gamma)]$-module on one generator.

Proof.  We define $U_n$ and $W_n$ to be the degree $n$ pieces of $U$
and $W$.  Let $w_\mu$ be the basis of elements
$\gamma(1)^{j_1}\gamma_2^{j_2}\cdots e^\beta$ for $W_n\otimes \Q$
parameterized by partitions $\mu=1^{j_1}2^{j_2}\cdots$ of $n$. The
$\Z$ module $W_n'$ spanned by the $w_\lambda$'s is not $W_n$ but has
index $\prod_{|\mu|=n} F(\mu)$ in it.  We let the elements
$L_\lambda=L_1^{i_1}L_2^{i_2}\cdots$ be the basis for the space
$U_n\otimes \Q$ indexed by partitions $\lambda
=1^{i_1}2^{i_2}\cdots$. The $\Z$ module $U_n'$ spanned by the
$L_\lambda$'s is not $U_n$ but has index $\prod_{|\lambda|=n}
F(\lambda)$ in it by lemma 6.2.  We define $m_{\lambda,\mu}$ for
$|\lambda|=|\mu|$ by $L_\lambda(w_\mu)=m_{\lambda,\mu}e^\beta$.  We
will show that the determinant of the $p(n)$ by $p(n)$ matrix
$(m_{\lambda,\mu})$ is
$$
\prod_{|\lambda|=n} (\beta,\gamma)^{l(\lambda)}P(\lambda)F(\lambda)
$$
where $l(\lambda)$ is the number of elements of the partition $\lambda$.
We order the partitions by $\lambda> \mu $ if 
$\lambda$ is the partition $\lambda_1+\lambda_2+\cdots$ with 
$\lambda_1\ge \lambda_2\ge\cdots$, 
$\mu$ is the partition $\mu_1+\mu_2+\cdots$ with 
$\mu_1\ge \mu_2\ge\cdots$, 
and $\lambda_1=\mu_1,\ldots,\lambda_{k-1}=\mu_{k-1}$, $\lambda_k>\mu_k$
for some $k$.
Then the matrix entry $m_{\lambda,\mu}$ is 0 if $\lambda>\mu$,
and is equal to
$$P(\lambda)F(\Lambda)(\beta,\gamma)^{l(\lambda)}$$
if $\lambda=\mu$. We can see this by repeatedly using the relation
$$\eqalign{
&\cdots L_{i_{2}} L_{i_1}
\gamma(j_1)\gamma(j_{2})\cdots e^\beta\cr
=&
\cases{
(\hbox{number of $j$'s equal to $j_1$}) (\beta,\gamma) j_1
\cdots L_{i_{2}}
\gamma(i_{2})\cdots e^\beta
&if $i_1=j_1$\cr
0&if $i_1>j_1$\cr}
\cr
}$$ 
for $j_1\ge j_2\ge \cdots$, which in turn follows from the
identities $[L_i,\gamma(j)]=j\gamma(j-i)$,
$L_i(e^\beta)=0=\gamma(-i)(e^\beta)$ for $i>0$, and
$\gamma_0e^\beta=(\beta,\gamma)$ As the matrix $(m_{\lambda,\mu})$ is
triangular its determinant is given by the product of the diagonal
entries
$m_{\lambda,\lambda}=P(\lambda)F(\Lambda)(\beta,\gamma)^{l(\lambda)}$.

Now we work out the index of $U^+(Witt_{>0})(e^{\beta^*})_n$ in
$W^*_n$ where $e^{\beta*} $ is the basis element of $W_0^*$ dual to
$e^\beta\in W_0$.  This index is equal to
$${
(\hbox{Index of $U'_n$ in $W'_n$})
\over
(\hbox{Index of $U'_n$ in $U_n$})(\hbox{Index of $W'_n$ in $W_n$})
}
$$
The numerator of this expression is equal to the determinant
of the matrix $(m_{\lambda,\mu})$, and we calculated this earlier.
Substituting in the known values we find the index
of $U^+(Witt_{>0})(e^{\beta*})_n$ in $W^*_n$ is
$${
\prod_{|\lambda|=n} P(\lambda)F(\lambda)(\beta,\gamma)^{l(\lambda)}
\over
(\prod_{|\lambda|=n} F(\lambda))(\prod_{|\mu|=n} F(\mu))
}$$
Applying lemma 6.1 we see that this is equal to
$\prod_{|\lambda|=n}(\beta,\gamma)^{l(\lambda)}$. This is a unit
in $\Z[1/(\beta,\gamma)]$, so over the ring $\Z[1/(\beta,\gamma)]$
the map from $U^+(Witt_{>0})$ to $W^*$ is an isomorphism.
This proves lemma 6.3.

Fix a norm 0 vector $\gamma\in L$. 
We recall that the transverse space if the subspace
of elements $v\in V$ such that $L_i(v)=0$ for $i>0$, 
$L_0(v)=v$, and $\gamma(i)(v)=0 $ for $i<0$. 
It is easy to check that the transverse space $T_\beta$ of degree 
$\beta$ with $(\beta,\gamma)\ne 0$ is positive definite, and the no-ghost
theorem [G-T] works by showing that the natural map from 
$T_\beta$ to the space of physical states (modulo null vectors)
of degree $\beta$ is an isomorphism. 

\proclaim Lemma 6.4. The transverse space $T_\beta$ of degree
$\beta\in II_{1,1}$ has determinant dividing a power
of $(\beta,\gamma)$.

Proof. For any vector $v$ of $V$ there is a $U^+(Witt_{>0})$ equivariant
map from $W^*$ to $V$ taking $1$ to $v$ over $\Z[1/(\beta,\gamma)]$ by
lemma 6.3. Dualizing we see that this means there is a vector in
$W\otimes V$ of the form $v\otimes 1+$ (terms involving some
$\gamma(i)$) which is fixed by $U^+(Witt_{>0})$ and hence also fixed
by $Witt_{>0}$.  This gives a map from $V$ to the transverse space
which is an isomorphism over $\Z[1/(\beta,\gamma)]$. This isomorphism
preserves the bilinear form because $\gamma$ has norm 0 so all the
terms involving $\gamma(i)$ have zero inner product with all other
terms. Hence the transverse space is self dual over
$\Z[1/(\gamma,\beta)]$.  This proves lemma 6.4.

The next theorem is an extension of the no-ghost theorem from rational
vector spaces to modules over $\Z$. Recall that the usual no-ghost
theorem [G-T] says that if $\alpha\in II_{25,1}$ is nonzero then the
space of physical states (over \Q) of degree $\alpha$ is positive
definite and spanned by the transverse space, which has dimension
$p_{24}(1-\alpha^2/2)$. This describes the space of physical states as
a rational vector space, but it also has a natural lattice inside it and
we can ask about the structure of this lattice.

\proclaim Theorem 6.5.~(The no-ghost theorem over \Z). Suppose $\alpha$
is a nonzero vector of $II_{25,1}$ which is
$n$ times a primitive vector. Then the discriminant of the space of
physical states of degree $\alpha$ divides a power of $n$.
In particular the space of physical states is a self dual lattice
if $\beta$ is primitive.

Proof. For every prime $p$ coprime to $n$ we can find a norm zero
vector $\gamma$ with $(\gamma,\beta)$ coprime to $p$, because the norm
0 vectors of Leech type span $II_{25,1}$. By lemma 6.4 this implies
that the discriminant of the space of physical states of degree
$\beta$ is coprime to $p$. Hence this discriminant divides a power of
$n$.  This proves theorem 6.5.

The proof of this theorem implies that all spaces of physical states
of primitive vectors of the same norm are isomorphic over all $p$-adic
fields, in other words in the same genus. It is certainly not always
true that they are isomorphic over the integers. For example, for a
norm 0 vector the space of physical states is isomorphic to the
corresponding Niemeier lattice, and not all Niemeier lattices are
isomorphic.

\proclaim 7.~An application to modular moonshine.

In this section $\m$ will stand for the fake monster Lie algebra
over $\Z[1/2]$ (see [B-R], [B98]) rather than the fake monster Lie algebra.

We write $\Z_p$ for the ring of $p$-adic numbers. 
The paper [B98] showed that Ryba's modular moonshine conjectures
[R] for primes $p\ge 13$ were true provided the following assumption
was true:
\proclaim Assumption. If $m<p$ then the degree $(m,n)$ piece
of $\m\otimes \Z_p$ is self dual under
the natural bilinear form and isomorphic to $V_{mn}\otimes \Z_p$
as a $\Z_p$ module acted on by the monster.

In this section we will prove this assumption, thus completing
the proof of the modular moonshine conjectures for $p\ge 13$.
In fact we will prove the following slightly stronger theorem:

\proclaim Theorem 7.1. If $m$ and $n$ are not both divisible by $p$ then
the degree $(m,n)$ piece
of $\m\otimes \Z_p$ is self dual under
the natural bilinear form and isomorphic to $V_{mn}\otimes \Z_p$
as a $\Z_p$ module acted on by the monster.

Proof. If one of $m$ or $n$ is not divisible by $p$ then we can find a
norm 0 vector $\gamma$ with $(\gamma,\beta)$ coprime to $p$. Theorem
7.1 then follows from the integral no-ghost theorem 6.5, because the
discriminant of $\m_{(m,n)}$ is not divisible by $p$ and is therefore
a unit in $\Z_p$.  Note that the isomorphism given by theorem 6.5
preserves the action of the monster. This proves theorem 7.1.

\proclaim 8.~Open problems.

\item{1} The integral form for the fake monster Lie algebra
constructed in section 4 is probably not the best possible one.  Is
there an integral form which not only has a structural basis but also
has the property that the bilinear form on the primitive elements is
self dual over $\Z$?
\item{2} 
Can the integral no-ghost theorem 6.5 be extended to show that the
root spaces of non-primitive vectors are self dual, and in particular
is the bilinear form on $\m$ self dual? We have shown in theorem 6.5
that it is self dual on the root spaces of roots that are either
primitive or of norm 0.  This question may be irrelevant, because it
is possible to construct an integral form of the fake monster Lie
algebra with a self dual bilinear form using the BRST complex.  Over
the complex numbers this construction is described in [F-G-Z] and
[L-Z]. The construction can be done over the integers, and the BRST
complex turns out to be a vertex superalgebra isomorphic to the vertex
superalgebra of the odd self dual lattice $I_{26,1}$.  The fake
monster Lie algebra is given by the homology groups of the self
adjoint BRST operator $Q$ modulo torsion, and the homology modulo
torsion of a complex with a self dual bilinear form is automatically
self dual under the induced form. However it is not clear that this self dual 
integral form can be extended to a smooth Hopf algebra over $\Z$. 
\item{3} Which elements of $\m_\alpha$ are liftable, and in particular
are all elements liftable?
\item{4}
We can also ask all the questions above about the monster Lie algebra
rather than the fake monster Lie algebra. In this case much less is
known than for the fake monster Lie algebra. In particular I do not
even know of an integral self dual form of the monster vertex algebra
or the monster Lie algebra (but see [B-R] for some speculation about
this).
\item{5} 
The arguments in section 4 of [B98] proving the modular moonshine
conjectures for $p\ge 13$ are rather messy and computational, mainly
because of the lack of a good description of what happens at the root
spaces of roots $\alpha$ divisible by $p$. Is is possible to clean up
this argument using the (conjectured) existence of a nice Hopf algebra
for the monster Lie algebra?

\proclaim References.

\item{[A]} E. Abe, Hopf algebras.
Cambridge Tracts in Mathematics, 74. Cambridge University Press,
Cambridge-New York, 1980. xii+284 pp.
ISBN: 0-521-22240-0
\item{[B86]}{R. E. Borcherds, Vertex algebras, Kac-Moody algebras,
and the monster. Proc. Natl. Acad. Sci. USA. Vol. 83 (1986) 3068--3071.}
\item{[B90]} R. E. Borcherds, The monster Lie algebra,
Adv. Math. Vol. 83, No. 1, Sept. 1990.
\item{[B98]} R. E. Borcherds, Modular Moonshine III, to appear in
Duke Math Journal.
\item{[B-R]}{R. E. Borcherds, A. J. E. Ryba,
Modular Moonshine II, Duke Math Journal Vol. 83 No. 2, 435-459, 1996. }
\item{[D]} J. Dieudonn\'e, Groupes de Lie et hyperalg\`ebres de Lie
sur un corps de caract\'eristique $p>0$ (V),
Bull. Soc. Math. France 84 1956 207--239
reprinted in
J. Dieudonn\'e, ``Choix d'\oe uvres math\'ematiques. Tome II.''
Hermann, Paris, 1981, ISBN:\quad 2-7056-5923-4, p. 600--632.
\item{[D72]} E. J. Ditters, Curves and formal (co)groups,
Inv. Math. 17 (1972) 1-20.
\item{[F-G-Z]} I. Frenkel, H. Garland, G. J. Zuckerman,
Semi-infinite cohomology and string theory. Proc.
Nat. Acad. Sci. U.S.A. 83 (1986), no. 22, 8442--8446.
\item{[G-T]} P. Goddard and C. B. Thorn, Compatibility of the dual
Pomeron with unitarity and the absence of ghosts in the dual resonance
model, Phys. Lett., B 40, No. 2 (1972), 235-238.
\item{[H]} M. Hazewinkel, ``Formal groups and applications''.
Pure and Applied Mathematics, 78.
Academic Press, Inc. [Harcourt Brace Jovanovich, Publishers],
New York-London, 1978. ISBN: 0-12-335150-2
\item{[K]} B. Kostant, Groups over $Z$.
1966 Algebraic Groups and Discontinuous Subgroups (Proc.
Symposium. Pure Math., Boulder, Colorado., 1965) pp. 90--98
Amer. Math. Soc., Providence, R.I.
\item{[L-Z]} B. Lian, G. J. Zuckerman, Moonshine cohomology.
Moonshine and vertex operator
algebra (Kyoto, 1994).
$\hbox{S\=urikaisekikenky\=usho}$
$\hbox{K\=oky\=uroku}$ No. 904 (1995), 87--115.
\item{[M]} I. G. Macdonald, Symmetric functions and Hall polynomials.
Second edition. Oxford Mathematical Monographs. Oxford Science Publications.
The Clarendon Press, Oxford
University Press, New York, 1995. x+475 pp. ISBN: 0-19-853489-2
\item{[R]}{A. J. E. Ryba, Modular Moonshine?,
In ``Moonshine, the Monster, and related topics'',
edited by Chongying Dong and Geoffrey Mason.
Contemporary Mathematics, 193. American Mathematical Society,
Providence, RI, 1996. 307-336. }
\item{[Sh]} P. B. Shay, An obstruction theory for smooth formal group
structure,
Preprint, Hunter college, CUNY.
\bye